\documentclass{acm_proc_article-sp}

\newtheorem{lemma}{Lemma}
\newtheorem{theorem}[lemma]{Theorem}
\newtheorem{corollary}[lemma]{Corollary}
\newtheorem{proposition}[lemma]{Proposition}
\newtheorem{definition}{Definition}

\def\bt{\begin{theorem}\rm}
\def\et{\end{theorem}}
\def\bc{\begin{corollary}\rm}
\def\ec{\end{corollary}}
\def\bl{\begin{lemma}\rm}
\def\el{\end{lemma}}
\def\bp{\begin{proposition}\rm}
\def\ep{\end{proposition}}
\def\bd{\begin{definition}\rm}
\def\ed{\end{definition}}

\def\bd{\begin{description}}
\def\ed{\end{description}}
\def\bu{\begin{enumerate}}
\def\eu{\end{enumerate}}
\def\bi{\begin{itemize}}
\def\ei{\end{itemize}}

\newbox\bigstrutbox
\setbox\bigstrutbox=\hbox{\vrule height12.5pt depth5pt width0pt}
\def\bigstrut{\relax\ifmmode\copy\bigstrutbox\else\unhcopy\bigstrutbox\fi}

\newbox\Bigstrutbox
\setbox\Bigstrutbox=\hbox{\vrule height17.5pt depth5pt width0pt}
\def\Bigstrut{\relax\ifmmode\copy\Bigstrutbox\else\unhcopy\Bigstrutbox\fi}

\def\ba{\begin{array}}
\def\ea{\end{array}}
\def\be{\begin{equation}}
\def\ee{\end{equation}}
\def\e{{\bf e}}

\def\P{{\bf P}}
\def\8{{\bf 8}}
\def\9{{\bf 9}}
\def\0{{\bf 0}}
\def\1{{\bf 1}}
\def\2{{\bf 2}}
\def\3{{\bf 3}}
\def\4{{\bf 4}}
\def\5{{\bf 5}}
\def\6{{\bf 6}}
\def\7{{\bf 7}}

\def\arrow{\overrightarrow}

\def\ds{\displaystyle}

\setpagenumber{1}

\begin{document}

\title{A Recipe for Symbolic Geometric Computing: Long Geometric Product,
BREEFS and Clifford Factorization}

\numberofauthors{1}

\author{\alignauthor Hongbo Li\\ \vskip .1cm
       \affaddr{Mathematics Mechanization Research Center}\\
       \affaddr{AMSS, Chinese Academy of Sciences}\\
       \affaddr{Beijing 100080, China}\\
       \email{hli@mmrc.iss.ac.cn}
}

%\date{}
\maketitle

\begin{abstract} \vskip .1cm

In symbolic computing, a major bottleneck is middle expression swell.
Symbolic geometric computing based on invariant algebras can alleviate this difficulty. For example,
the size of projective geometric computing based on bracket algebra can often be restrained to
two terms, using final polynomials, area method, Cayley expansion, etc. This is the
``binomial" feature of projective geometric computing in the language of bracket algebra.

In this paper we report a stunning discovery in Euclidean geometric computing: the
term preservation phenomenon. Input an expression in the language of Null Bracket
Algebra (NBA), by the recipe we are to propose in this paper,
the computing procedure can often be controlled
to within the same number of terms as the input, through to the end. In particular,
the conclusions of
most Euclidean geometric theorems can be expressed by monomials
in NBA, and the expression size in the proving procedure can often be controlled to within one term!
Euclidean geometric computing can now be
announced as having a ``monomial" feature in the language of NBA.

The recipe is composed of three parts: use long geometric product to
represent and compute multiplicatively, use ``BREEFS" to control the expression size locally,
and use Clifford factorization for term reduction and transition from algebra to geometry.

By the time this paper is being written, the recipe has been tested by $70+$ examples from \cite{chou},
among which $30+$ have monomial proofs. Among those outside the scope, the famous
Miquel's five-circle theorem \cite{chou2}, whose analytic proof is straightforward but very difficult
symbolic computing, is discovered to have a 3-termed elegant proof with the recipe.

{\bf ACM Computing Classification:}\
I.1.1 [Symbolic and Algebraic Manipulation]: Expressions and Their Representation;
G.4 [Mathematical Software]: Efficiency.

{\bf General Terms:}\ Theory; algorithm.


\end{abstract}

\section{Introduction}
\vskip .3cm

Using geometric invariants in symbolic geometric computing has been an active research subject
in symbolic and algebraic computation. Apart from the benefit of better geometric
interpretability when compared with coordinates \cite{wang}, geometric invariants have a salient feature of
reducing the size of symbolic manipulation. This is particularly valuable because a major difficulty
in symbolic computing is middle expression swell.

In projective incidence geometry, the proofs of many
theorems by the method of biquadratic final polynomials \cite{gebert} can be so elegant
that only bracket binomials occur in the whole procedure. The area method
\cite{chou2} also shares this feature in its ratio-formed proofs of many theorems, i.e.,
the numerators and denominators are monomials of areas. The advantages of both methods are assimilated
into the Cayley expansion theory developed in \cite{li-wu}, by which this feature of projective
incidence geometry is extended to projective conic geometry.

In Euclidean geometry, with the introduction of inner products, syzygies among basic invariants
become much more complicated. Using distances for geometric theorem proving was proposed in \cite{havel},
and further developed in \cite{mourrain} by including the inner products of extensors. In \cite{li04}
the covariant algebra of Euclidean geometry, the so-called conformal geometric algebra (CGA), and its
invariant subalgebra, the so-called null bracket algebra, was shown to provide a nice algebraic setting for
Euclidean geometric theorem proving. Unfortunately, except for some sporadic special cases,
none of these methods has ever shown any binomial feature. As a benchmark problem,
Miquel's 5-circle theorem \cite{chou2}, whose analytic proof is straightforward but very difficult
symbolic computing,
was found a proof of 14 terms in 2001 \cite{li05}. The proof was full of
pairwise term reductions based on complicated syzygies of advanced invariants.

Under such background it comes as appalling as can be an observation that symbolic
computing and theorem proving
in Euclidean geometry have a ``monomial" feature, or more generally, a
term preservation feature in the language of NBA, if the recipe in this paper is used.

By the time this paper is being written, $70+$ Euclidean geometric theorems in \cite{chou} have been
tested, and $40+$ have the term preservation feature. In particular, $30+$ theorems have their
conclusions represented by monomials in NBA, and are kept as monomials till the end of the proof.
Concrete examples include all the examples published in \cite{li04}. Using the method in that paper
only one example preserves its number of terms in the proof.
Using our recipe in this paper ALL examples preserve their number of terms, thus it is impossible to
find any analytic proof that is more elegant.

\vskip -.34cm
\begin{table}[htbp]
\centering
\caption{Comparison of proofs with [8]}
\vskip .15cm
 \begin{tabular}{|c|c|c|c|}
\hline
Example in [8] & Conclusion & Proof in [8] & New proof\\ \hline
No. 1 & 1 term & 3 terms & 1 term\\ \hline
No. 2 & 2 terms & 2 terms & 2 terms \\ \hline
No. 3 & 1 term & 3 terms & 1 term \\ \hline
No. 4 & 2 terms & 3 terms & 2 terms \\ \hline
No. 5 & 1 term & 4 terms & 1 term \\ \hline
No. 6 & 1 term & 4 terms & 1 term \\ \hline
\end{tabular}
\label{tab}
\end{table}

By preserving the number of terms the computing burden is transmitted from addition to noncommutative
multiplication.
One may shake head as to any possible simplification by algebraic manipulation of multiplication in place
of addition. Well, in an invariant symbolic
system there are syzygies among basic elements.
In manipulating such elements, multiplication
preserves geometry while addition breaks it up.
Symmetries in multiplication provide the most economical way of avoiding or employing syzygies.
It is easy to change multiplication to addition: just recall how coordinates
are introduced. Generally it is very difficult to change addition
to multiplication: just recall Cayley factorization \cite{white} in projective geometry.

Our recipe for symbolic computing in Euclidean geometry is: (1)
employ multiplication, or more accurately, the {\it geometric product} in {\it Geometric Algebra},
from the representation of geometric objects on, (2) preserve the multiplication through
subsequent algebraic manipulation using the principle ``BREEFS"
\cite{li-wu}, and (3)
replace addition by multiplication using Clifford factorization -- the Euclidean version of
Cayley factorization. In (1) we need to invent two new devices
for the representation by multiplication, called {\it nullifying operator} and {\it reduced meet product}.
In (2) we need to adapt
the previous global invariant bracket-oriented principle to a
shift invariant {\it neighborhood principle}. In (3)
we need a device to explore rational Clifford expansions systematically
-- {\it pseudodivision} in NBA. These novelties will be introduced
in Sections 4 and 5, with various illustrations.

Geometric Algebra \cite{hestenes} is a version of Clifford algebra favoring the universal usage of its
multiplication, the {\it geometric product}, instead of addition. Hestenes' vision of Geometric Algebra
in place of the more commonly used Clifford algebra in matrix or hypercomplex numbers form,
is fully justified by our theorem proving practice: replacing
addition by multiplication and prolonging the multiplication (``{\it long geometric product}"),
are the simplest means of avoiding syzygies because the geometric product already
incorporates various syzygies of inner products and determinants into its structural
symmetry.

The Geometric Algebras developed for the conformal model of Euclidean geometry,
CGA \cite{li01} and NBA \cite{li02},
will be introduced in Section 2 from the implementation point of view. The geometry of
long geometric product in NBA will be explained in Section 3. In
the end of this paper, a 3-termed analytic proof will be provided for
the benchmark problem, Miquel's 5-circle theorem.

\section{Conformal geometric algebra \\
and null bracket algebra}
\vskip .3cm

In \cite{li01}, \cite{li02}, \cite{li04} there have been detailed introductions of CGA and NBA.
In this paper we concentrate on the case of 2D geometry only. We always use
a boldfaced number or letter to denote a vector.

In 4D Minkowski space we fix a null vector $\e$. A {\it null vector} is a nonzero vector whose inner product
with itself is zero. We call $\e$ the {\it point at infinity} of the Euclidean plane. Any null vector
linear independent of $\e$ is a {\it point} in the plane. Two null vectors represent the same point
if and only if they differ only by scale. This representation of the Euclidean plane is {\it conformal} but not
isometric. To obtain an isometric model we only need to fix the inner product of any point with the point
at infinity, e.g. to $-1$, as any two linear independent null vectors have nonzero inner product.

To describe and analyze Euclidean geometry with the conformal model we need a suitable algebraic language.
The symbolic version of Clifford algebra in \cite{hestenes} is an ideal tool in that
it prefers the usage of multiplication to addition. The multiplication, called geometric product,
conglomerates all geometric relations within itself and is geometrically meaningful. The other
two versions, the matrix version and the hypercomplex numbers version, emphasize the linear nature of the
algebra, i.e., care more for addition than for multiplication. In symbolic form, more addition
leads to more algebra, and more multiplication preserves more geometry. This justifies the gist ``geometric"
in Hestenes' Geometric Algebra.

{\it Geometric Algebra} is the unique algebra generated from an $n$D
inner product space by an associative product, called {\it geometric product},
satisfying the generating relation that the geometric product of any vector with itself is the inner product.
The geometric product is always denoted by juxtaposition.

Geometric Algebra is graded, the {\it grade} ranges from 0 to $n$. For an element $A$ in this algebra,
the {\it $i$-graded part} is denoted by $\langle A \rangle_i$. When $i=0$ and 1 it is the
{\it scalar} and {\it vector part}
respectively. Elements of grade $i$ form a subspace of dimension $C_n^i$. In particular when $i=n$, the
{\it $n$-graded subspace} is 1D. Fix a nonzero element $I_n$ in this space, then for any other $n$-graded
element $A_n$, the coordinate of $A_n$ with respect to the basis $I_n$ is the {\it bracket} of $A_n$:
\be
[A_n]=\frac{A_n}{I_n}=A_nI_n^{-1}=I_n^{-1}A_n.
\ee
The geometric product of an element $A$ with $I_n^{-1}$ is called the {\it dual} of $A$: \
$A^\sim=AI_n^{-1}$. In particular $A_n^\sim=[A_n]$.

The geometric product of two vectors is composed of two parts, the 0-graded part and the
2-graded part. They are respectively the {\it inner product} and {\it outer product}, denoted by
dot and wedge. The outer product is just the exterior product in Grassmann's exterior algebra.
\be
\1\2+\2\1=2\,\langle \1\2\rangle=2\,(\1\cdot \2),\ \ \
2\,(\1\wedge \2)=\1\2-\2\1.\label{innout}
\ee
For three vectors,
\be
2\,\langle \1\2\3 \rangle_1=\1\2\3-\3\2\1,\ \ \
2\,(\1\wedge \2\wedge \3)=\1\2\3-\3\2\1.\label{innout3}
\ee

\vskip -.2cm
{\it Conformal geometric algebra} (CGA) \cite{li01} is the Geometric Algebra established upon the conformal model.
In CGA for 2D geometry, we can pick out two scalar-valued functions generated by
a sequence of $2k$ null vectors $\1, \2, \cdots$, ${\bf 2k}$: the 0-graded part
\be
\langle \1\2\cdots ({\bf 2k})\rangle = \langle \1\2\cdots ({\bf 2k})\rangle_0,
\ee
called the {\it angular bracket},
and the bracket of the 4-graded part for $k>1$:
\be
[\1\2\cdots ({\bf 2k})] = [\langle \1\2\cdots ({\bf 2k})\rangle_4],
\ee
called the {\it square bracket}.
The two kinds of brackets generate a ring called {\it null bracket algebra} (NBA) \cite{li02}, \cite{li04}.

NBA should be implemented by realizing
the following basic properties:

(1) Multilinearity of geometric product, multilinearity and
(anti-)symmetry of inner product and outer product, (\ref{innout}) and (\ref{innout3})
from left to right.\\
\\
(2) Null symmetry: for null vector $\1$,
\be
\1\2\3\1=-\1\3\2\1=\1(\2\wedge \3)\1, \ \ \ \1\1=\1\cdot \1=0. \label{cga:nullsymmetry}
\ee
(3) Shift and reversion symmetry: for vectors $\bf 1,2,\cdots, k$,
\be\ba{lll}
\,[{\bf 12\cdots k}] &=& -[{\bf 2\cdots k1}]\, = -[{\bf k1\cdots (k-1)}],\\
\langle {\bf 12\cdots k}\rangle
&=&\  \ \langle{\bf 2\cdots k1}\rangle =\ \  \langle{\bf k1\cdots (k-1)}\rangle,\\
\,{[}{\bf 12\cdots k}] &=& \ \ [{\bf k(k-1)\cdots 21}],\\
\langle {\bf 12\cdots k}\rangle &=& \ \, \langle {\bf k(k-1)\cdots 21}\rangle.
\ea
\label{cga:shift}
\ee
(4) Dual symmetry: for element $A$ and $r$-graded element $B_r$,
\be\ba{l}
\1(\1\wedge \2\wedge \3)^\sim=-(\1\wedge \2\wedge \3)^\sim\1,\\
 \langle A^\sim \rangle=[A], \ \ \
[A^\sim]=-\langle A\rangle, \\
A^\sim B_r=(-1)^r (A B_r)^\sim, \ \ \ {}^{\sim\sim}=-1.
\ea
\label{dualsymmetry}
\ee
(5) Points $\1,\2,\3$ are collinear:
$[\e\1\2\3]=0;$
points $\1,\2,\3,\4$ are cocircular:
$[\1\2\3\4]=0.$
Lines $\1\2, \1'\2'$ are parallel:\\
$[\e\1\2\e\1'\2']=0$; they are perpendicular:
$\langle \e\1\2\e\1'\2'\rangle=0$.\\
\\
(6) Contraction by Grassmann-Pl\"ucker syzygy:
for vectors $\1,\2,\3,\4,\5,\1',\2',\3'$,
\be\hskip -.1cm
\ba{r}
[\1\2\3\4][\5\1'\2'\3']-[\1\2\3\5][\4\1'\2'\3']+[\1\2\4\5][\3\1'\2'\3']\\
=[\1\3\4\5][\2\1'\2'\3']
-[\2\3\4\5][\1\1'\2'\3'].
\ea
\label{cga:gp}
\ee
(7) Contraction by inner-product bracket syzygy:
for vectors $\1, \2$, $\3, \4, \5, \1'$,
\be\ba{r}
\1'\cdot \1[\2\3\4\5]-\1'\cdot \2[\1\3\4\5]
+\1'\cdot \3[\1\2\4\5]\\
\hskip .5cm =\1'\cdot \4[\1\2\3\5]
-\1'\cdot \5[\1\2\3\4].
\ea
\ee
(8) Null expansion:
for null vector $\1$,
\be\ba{lll}
\1\2\1 &=& 2\,(\1\cdot \2)\1,\\
\1\2\3\1 &=& 2\,(\1\cdot \2)\3\1-2\,(\1\cdot \3)\2\1.
\ea
\label{cga:nullexpansion}
\ee
(9) Trigonometric
quartet expansion: for null vector $\1$,
\be\ba{rcl}
\ds\frac{1}{2}[\1\2\3\4\1\5\cdots {\bf r}]
&=& \langle \1\2\3\4 \rangle\,[\1\5\cdots {\bf r}]\\
&& \hskip -.27cm
+\,[\1\2\3\4]\,\langle \1\5\cdots {\bf r}\rangle, \\
\ds\frac{1}{2}\langle \1\2\3\4\1\5\cdots{\bf r} \rangle\bigstrut
&=& \langle \1\2\3\4 \rangle\, \langle\1\5\cdots {\bf r}\rangle\\
&&\hskip -.27cm
-\,[\1\2\3\4]\,[\1\5\cdots {\bf r}].
\ea
\label{quartet}
\ee
(10) Trigonometric
sextet expansion: for null vector $\1$,
\be
\ba{lcl}
\ds\frac{1}{2}[\1\2\3\1\4\5\1\6\cdots {\bf r}]
&=&\ [\1\2\3\1\4\5]\langle\1\6\cdots {\bf r}\rangle\\
&&\hskip -.191cm
+\langle \1\2\3\1\4\5 \rangle[\1\6\cdots {\bf r}], \\
\bigstrut
\ds\frac{1}{2}\langle \1\2\3\1\4\5\1\6\cdots {\bf r}\rangle
&=& \langle \1\2\3\1\4\5\rangle \langle\1\6\cdots {\bf r}\rangle\\
&&
\hskip -.22cm -[\1\2\3\1\4\5][\1\6\cdots {\bf r}].
\ea
\label{li}
\ee
(11) Rational sextet expansion: for null vectors $\2,\3,\5,\6$,
\be
%\hskip -.32cm
\ba{rcl}
\ds -\frac{1}{2}[\1\2\3\4\5\6][\2\3\5\6]
&=& \2\cdot \3[\1\2\5\6][\3\4\5\6]\\
&& \hskip -.23cm  +\5\cdot \6[\1\2\3\6][\2\3\4\5],\\
\ds\bigstrut -\frac{1}{2}\langle \1\2\3\4\5\6\rangle[\2\3\5\6]
&=& \2\cdot \3\,[\1\2\5\6]\langle \3\4\5\6\rangle\\
&&\hskip -.23cm -\5\cdot \6\langle \1\2\3\6\rangle[\2\3\4\5].
\ea
\label{cga:miquel4}
\ee
(12) Rational octet expansion: for null vectors $\1,\2,\3,\5$,
\be
\hskip -.2cm
\ba{lcl}
\ds \frac{1}{2}[\1\2\3\4\1\2\5\6][\1\2\3\5]
&=& \2\cdot \3[\1\2\5\6][\1\2\5\1\3\4]\\
&&\hfill -\2\cdot \5[\1\2\3\4][\1\2\3\1\5\6],\\
\ds \frac{1}{2}\langle \1\2\3\4\1\2\5\6 \rangle[\1\2\3\5]\bigstrut
&=& \2\cdot \3\langle \1\2\5\6 \rangle[\1\2\5\1\3\4]\\
&&\hfill\ -\2\cdot \5[\1\2\3\4]
\langle \1\2\3\1\5\6\rangle,\\
\ds\Bigstrut \frac{1}{2}[\1\3\2\4\1\5\2\6][\1\3\2\5] &=& (\1\cdot \5)(\2\cdot \5)[\1\3\2\4][\1\3\2\6]\\
&&\hskip -.23cm
+(\1\cdot \3)(\2\cdot \3)[\1\5\2\4][\1\5\2\6],\\
\ds\bigstrut  \frac{1}{2}\langle \1\3\2\4\1\5\2\6 \rangle [\1\3\2\5]
&=& (\1\cdot \5)(\2\cdot \5)[\1\3\2\4]\langle \1\3\2\6\rangle\\
&& \hskip -.23cm
-(\1\cdot \3)(\2\cdot \3)[\1\5\2\4]\langle \1\5\2\6\rangle.
\ea
\label{cga:miquelformula}
\ee
(13) Reverse of the expansions from (\ref{quartet}) to ({\ref{cga:miquelformula}):
from right to left. They are basic {\it Clifford factorizations}.

The derivation of the properties is easy. In Section \ref{sect:geom},
(\ref{quartet}) and (\ref{li}) are discussed, in Section \ref{sect:breefs},
 (\ref{cga:miquel4}) and (\ref{cga:miquelformula}) are analyzed.

\section{the geometry of long product}
\label{sect:geom}
\vskip .3cm

Prolonging the length of elements in inner products and brackets (i.e., determinants) is not
only a device of simplifying symbolic computing, but an indispensable means of representing basic
geometric relations. We take a look at a planar angle and its algebraic representation.

An oriented angle is just a 2D rotation. $\angle \1\2\3$ can be represented by three points:
the vertex $\2$, a point $\1$ on the initial ray, and a point $\3$ on the terminal ray without
requiring that $\1,\3$ be equidistant from $\2$.
Without resorting to inequalities, any rational description of angle $\angle \1\2\3$ by points $\1,\2,\3$
is accurate only up to $k\pi$.
The equivalent classes of oriented planar angles modulo $\pi$ are called {\it full angles} \cite{chou2}.

Let $\angle \1\2\3,\angle \1'\2'\3'$ be two full angles. They are equal if and only if
$\tan\angle \1\2\3=\tan\angle \1'\2'\3'$. In the conformal model, the ratio of
$[\e\1\2\3]$ to $\langle \e\1\2\3 \rangle$ is exactly
$\tan\angle \1\2\3$. So $\angle \1\2\3=\angle \1'\2'\3'$ if and only if
\be
\frac{[\e\1\2\3]}{\langle \e\1\2\3\rangle}=
\frac{[\e\1'\2'\3']}{\langle \e\1'\2'\3'\rangle}, \label{full}
\ee
i.e.,
\be
\frac{1}{2}[\e\1\2\3\e\3'\2'\1']=
[\e\1\2\3]\langle \e\1'\2'\3'\rangle-\langle \e\1\2\3\rangle[\e\1'\2'\3']
=0.
\ee

Essentially, geometric product $\e\1\2\3$ represents full angle $\angle \1\2\3$, its sine and cosine are
respectively the square and angular brackets, its reverse angle is $\e\3\2\1$.
$\e\1\2\e\3\4$ represents full angle $\angle(\overrightarrow{\bf 12},\overrightarrow{\bf 34})$.
The sum of full angles $\e\1\2\3, \e\1'\2'\3'$ is their concatenation
$\e\1\2\3\e\1'\2'\3'$. These explain (\ref{quartet}) and (\ref{li}) as
expansions of the sines and cosines of angle sums.

Using null expansion (\ref{cga:nullexpansion}) and trigonometric expansions (\ref{quartet}) and (\ref{li}),
we easily obtain trigonometric explanations of all square and angular
brackets. For example if $\e$ does not occur in the following vectors, then
\[
\ba{lll}
\langle \1\2\3\4 \rangle
&=& \ds -\frac{d_{\bf 12}d_{\bf 23}d_{\bf 34}d_{\bf 41}}{2}\cos \angle(\arrow{\1\2\3}, \arrow{\1\3\4}),\\
{[}\1\2\3\4]
&=& \ds -\frac{d_{\bf 12}d_{\bf 23}d_{\bf 34}d_{\bf 41}}{2}\sin \angle(\arrow{\1\2\3}, \arrow{\1\3\4}),
\Bigstrut \\
\langle \1\2\3\4\5\6 \rangle\Bigstrut
&=& \ds -\frac{d_{\bf 12}d_{\bf 23}d_{\bf 34}d_{\bf 45}d_{\bf 56}d_{\bf 61}}{2}\cos
(\angle(\arrow{\1\2\3}, \arrow{\1\3\4})\\
&&\hfill +\angle(\arrow{\1\4\5}, \arrow{\1\5\6})),\\
{[}\1\2\3\4\5\6]
&=& \ds -\frac{ d_{\bf 12}d_{\bf 23}d_{\bf 34}d_{\bf 45}d_{\bf 56}d_{\bf 61}}{2}\sin
(\angle(\arrow{\1\2\3}, \arrow{\1\3\4})\bigstrut\\
&&\hfill +\angle(\arrow{\1\4\5}, \arrow{\1\5\6})),
\ea
\]
and for the general case,
\[
\ba{ll}
& \langle {\bf i}_1{\bf i}_2\cdots {\bf i}_{2l+2} \rangle\\
=& \ds -\frac{d_{{\bf i}_1{\bf i}_2} d_{{\bf i}_2{\bf i}_3}\cdots d_{{\bf i}_{2l+1}{\bf i}_{2l+2}}
d_{{\bf i}_{2l+2}{\bf i}_1}}{2}\cos
(\angle(\arrow{{\bf i}_1{\bf i}_2{\bf i}_3}, \arrow{{\bf i}_1{\bf i}_3{\bf i}_4})\\
& \hfill \bigstrut
+\angle(\arrow{{\bf i}_1{\bf i}_4{\bf i}_5}, \arrow{{\bf i}_1{\bf i}_5{\bf i}_6})+\cdots+
\angle(\arrow{{\bf i}_1{\bf i}_{2l}{\bf i}_{2l+1}}, \arrow{{\bf i}_1{\bf i}_{2l+1}{\bf i}_{2l+2}}));\\
\\
& {[}{\bf i}_1{\bf i}_2\cdots {\bf i}_{2l+2}]\\
= &\ds -\frac{d_{{\bf i}_1{\bf i}_2} d_{{\bf i}_2{\bf i}_3}\cdots d_{{\bf i}_{2l+1}{\bf i}_{2l+2}}
d_{{\bf i}_{2l+2}{\bf i}_1}}{2}\sin\Bigstrut
(\angle(\arrow{{\bf i}_1{\bf i}_2{\bf i}_3}, \arrow{{\bf i}_1{\bf i}_3{\bf i}_4})\\
& \hfill \bigstrut
+\angle(\arrow{{\bf i}_1{\bf i}_4{\bf i}_5}, \arrow{{\bf i}_1{\bf i}_5{\bf i}_6})+\cdots+
\angle(\arrow{{\bf i}_1{\bf i}_{2l}{\bf i}_{2l+1}}, \arrow{{\bf i}_1{\bf i}_{2l+1}{\bf i}_{2l+2}})).
\ea
\]
Here $d_{\1\2}$ is the Euclidean distance between points $\1,\2$;
$\arrow{\1\2\3}$ denotes the oriented circle through $\1,\2,\3$ sequentially, and
$\angle(\arrow{\1\2\3}, \arrow{\1\3\4})$ is the full angle from the tangent direction of
$\arrow{\1\2\3}$ to that of $\arrow{\1\3\4}$ at any point of their intersection.

If $\e$ occurs then the explanation is only slightly changed. For example,
\be\ba{lll}
\langle \e\1\2\3\4\5\rangle &=& -d_{\bf 12}d_{\bf 23}d_{\bf 34}d_{\bf 45}\cos (\angle{\1\2\3}+ \angle{\3\4\5}),\\
\,{[}\e\1\2\3\4\5] &=& -d_{\bf 12}d_{\bf 23}d_{\bf 34}d_{\bf 45}\sin (\angle{\1\2\3}+ \angle{\3\4\5}),
\ea
\ee
and the general case follows similarly. The power of long geometric product
comes from its intrinsic geometric nature.

\section{the power of long product}
\vskip .3cm

Below we present two novel devices in NBA, the {\it nullifying operator} and the
{\it reduced meet product}. They function as the bridge
between Grassmann-Cayley algebra and NBA, thus allowing to employ the full power of
Cayley expansion theory \cite{li-wu} within Euclidean geometry.

Let $\1$ be a null vector in 4D Minkowski space. For vector $\2$,
\be
N_\1(\2)=\frac{1}{2}\,\2\1\2
\ee
is the {\it nullification} of $\2$ with respect to $\1$. When $\2$ is null then
$N_\1(\2)=(\1\cdot \2)\2$ represents the same point $\2$. When $\2$ is not null, then
$N_\1(\2)$ represents the null vector other than $\1$ in the plane spanned by $\1, \2$
if the metric of the plane is
Minkowski, or just $\1$ if the metric is degenerate.

The {\it reduced meet product} of two elements $\2\wedge \3$ and $\2'\wedge \3'$ modulo
vector $\1$ is
\be\ba{lll}
(\2\wedge \3)\vee_\1 (\2'\wedge \3') &=& [\1\2\2'\3']\3-[\1\3\2'\3']\2 \\
&=& [\1\2\3\3']\2'-[\1\2\3\2']\3'.
\ea
\label{cga:meet}
\ee
The second equality is modulo $\1$, i.e., the two sides differ by $\lambda\,\1$ for a scale $\lambda$.
This product is a reduced form of the classical {\it meet product}
%``$\vee$"
of elements
$\1\wedge \2\wedge \3$ and $\1\wedge \2'\wedge \3'$:
\be
\hskip -.2cm
\ba{r}
(\1\wedge \2\wedge \3)\vee (\1\wedge \2'\wedge \3') = [\1\2\2'\3']\1\wedge \3-[\1\3\2'\3']\1\wedge \2 \\
= [\1\2\3\3']\1\wedge \2'-[\1\2\3\2']\1\wedge \3'
\ea\ee

{\it Proposition 1}. [Null duality]\ Let $\1'$ be a null vector, then
\be
\1'(\1\wedge \2\wedge \3)^\sim (\1\wedge \2'\wedge \3')^\sim \1'
=\1'\1\{(\2\wedge \3)\vee_\1(\2'\wedge \3')\}^\sim\1'.
\label{proposition1}
\ee

{\it Proof.} In Geometric Algebra we have the duality relation
\be
A^\sim\wedge B^\sim=(A\vee B)^\sim,
\ee
so
\[\ba{ll}
& \1'(\1\wedge \2\wedge \3)^\sim (\1\wedge \2'\wedge \3')^\sim \1'\\
=& \1'\{(\1\wedge \2\wedge \3)^\sim \wedge  (\1\wedge \2'\wedge \3')^\sim\} \1'\\
=& \1'\{(\1\wedge \2\wedge \3)\vee  (\1\wedge \2'\wedge \3')\}^\sim \1'\\
=& \1'\1\{(\2\wedge \3)\vee_\1(\2'\wedge \3')\}^\sim\1'.
\ea
\]

{\bf Example 1.} [See \cite{li04}, Example 5]\ If three circles having a point in common
intersect pairwise at three collinear points, their common point
are cocircular with their centers.

\begin{figure}[htbp]
\centerline{\epsfig{file=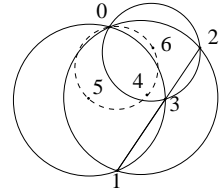, height=2.5cm}}
\caption{Example 1.}
\end{figure}

This is the most difficult example in \cite{li04}.
We use the same geometric scenario: remove
the collinearity constraint from the hypotheses, compute the
conclusion expression to see if the removed constraint
comes out as a factor.

Same as in \cite{li04}, during the computing all intermediate
factors are saved in a set for later analysis. They are marked
with under braces and are removed from subsequent steps.

Free points: $\bf 0, 1, 2$.
\\ \indent
Centers: $\4={\rm center}(\0\1\2),\ \5={\rm center}(\0\1\3)$,\ $\6={\rm center}(\0\2\3)$.
\\ \indent
Removed hypothesis: $[\e\1\2\3]=0$.
\\ \indent
Conclusion expression: $[{\bf 0456}]$.

The center of circle $\1\2\3$ is
\be
{\rm center}(\1\2\3)=N_\e((\1\wedge \2\wedge \3)^\sim).
\ee
Substitute the expressions of the three circle centers $\4,\5,\6$ into the conclusion, we get
\[
\hskip -.5cm
\ba{lrl}
& 2^3{[}\0\4\5\6]
=& [\0(\0\wedge \1\wedge \2)^\sim \e (\0\wedge \1\wedge \2)^\sim\bigstrut
(\0\wedge \1\wedge \3)^\sim \e \\
&& \hfill (\0\wedge \1\wedge \3)^\sim
(\0\wedge \2\wedge \3)^\sim \e (\0\wedge \2\wedge \3)^\sim ]\\
&=& [\0\e(\0\wedge \1\wedge \2)^\sim (\0\wedge \1\wedge \3)^\sim\e \bigstrut
(\0\wedge \1\wedge \3)^\sim\\
&&\hfill
(\0\wedge \2\wedge \3)^\sim \e (\0\wedge \2\wedge \3)^\sim
(\0\wedge \1\wedge \2)^\sim ]\\
&=& \bigstrut
[\0\e\0\{(\1\wedge \2)\vee_\0 (\1\wedge \3)\}^\sim\e\0
\{(\1\wedge \3)\vee_\0 \\
&&\hfill (\2\wedge \3)\}^\sim \e (\0\wedge \2\wedge \3)^\sim
(\0\wedge \1\wedge \2)^\sim ]\\
&=& \bigstrut\hskip -.2cm
\underbrace{-2\,(\e\cdot \0)[\0\1\2\3]^2}
[\0\1\e\0\3\e (\0\wedge \2\wedge \3)
(\0\wedge \1\wedge \2)]\\
&=& \, 2^{-2}\, [\0\1\e\0\3\e\0\2\3\0\1\2]\bigstrut\\
&=& \hskip -.2cm-2^{-2}\, [\0\1\e\0\e\3\0\3\2\0\2\1]\bigstrut\\
&=& \underbrace{2^2\,(\e\cdot \0)(\0\cdot \1)(\0\cdot \2)(\0\cdot \3)}[\e\1\2\3].\bigstrut
\ea
\label{cga:power}
\]
{\it Explanation of the computing}:
%in (\ref{cga:power})

Line 1: substitution.\\
Line 2: change of order by symmetries
(\ref{dualsymmetry}) and (\ref{cga:shift}).\\
Line 3: apply (\ref{proposition1}) to the first two pairs of meet products.\\
Line 4: expand meet products, make null expansion (\ref{cga:nullexpansion}).\\
Line 5: apply (\ref{innout3}). Only one term is generated.\\
Line 6: change of order between neighboring $\0$'s by
(\ref{cga:nullsymmetry}). \\
Line 7: null expansion.

Now consider the most typical geometric construction:
the intersection of two circles. Let $\1\2\3$ and $\1\2'\3'$ be two circles represented by
circumpoints. Their {\it point of intersection}, denoted by
$\1\2\3\cap \1\2'\3'$, refers to the point of intersection other than $\1$, or
in the tangent case, tangent point $\1$ itself. When $\1=\e$ it is the intersection of
lines $\2\3, \2'\3'$; when $\2=\e$ it is the intersection of line $\1\3$ and circle $\1\2'\3'$.

In
\cite{li02} the intersection was expressed as a linear combination of three vectors:
either $\1,\2,\3$, or $\1,\2',\3'$.
In this paper we propose the following representation:
\be
\1\2\3\cap \1\2'\3'=N_\1((\2\wedge \3)\vee_\1 (\2'\wedge \3')).
\label{cga:renovation}
\ee
It is easily verified that (\ref{cga:renovation})
equals the expressions in \cite{li02}.

In (\ref{cga:renovation}) the reduced meet product has two ways of expansion, either by
separating $\2,\3$ as in the first line of (\ref{cga:meet}), or by
separating $\2',\3'$. If we compute the geometric product
\be\hskip -.1cm
\4(\1\2\3\cap \1\2'\3')\5=\frac{1}{2}\4\{(\2\wedge \3)\vee_\1 (\2'\wedge \3')\}
\1\{(\2\wedge \3)\vee_\1 (\2'\wedge \3')\}\5, \label{cga:order}
\ee
then previously we simply expanded the meet products in the same way and multiplied them with $\1$.
Since the geometric product is associative, not only can we expand the meet products in different ways, but
we can freely change the order of the four pairwise geometric products in (\ref{cga:order}).
Furthermore, if the intersection has more than one construction, e.g. it is where three lines meet,
we can even change the representations by different pairs of lines for each meet product.
By prolonging the length of the geometric product, we gain a lot of freedom for the realization of
``BREEFS".

\section{BREEFS}
\label{sect:breefs}
\vskip .3cm

BREEFS -- ``{\bf B}racket-oriented {\bf R}epresentation, {\bf E}limination and {\bf E}xpansion for
{\bf F}actored and {\bf S}hortest result", was first proposed in \cite{li-wu}
to control the expression size in bracket algebra. In \cite{li04} it was
extended to null inner-product bracket algebra, the subalgebra of NBA generated by
angular brackets of length 2 and square brackets of length 4. In this section we extend
it to long geometric products. The best way to explain this principle is through some working examples.

{\bf Example 2.}
In the plane two circles intersect at points $\1,\1'$.
Two secant lines through them intersect the circles at points
$\2,\3$ and $\2',\3'$ respectively. Then $\2\2'//\3\3'$.

\begin{figure}[htbp]
\centering \epsfig{file=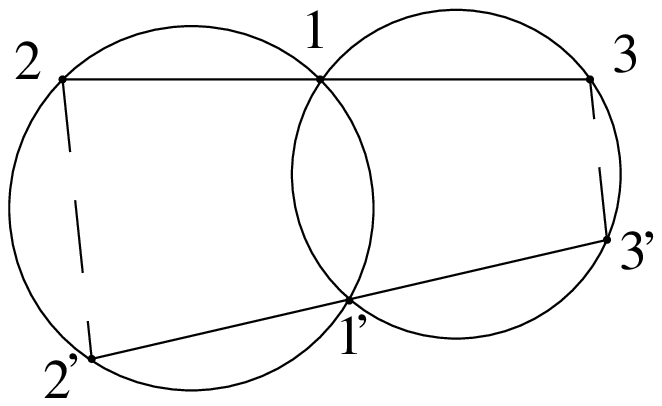, height=2cm}(a)\ \ \
\centering \epsfig{file=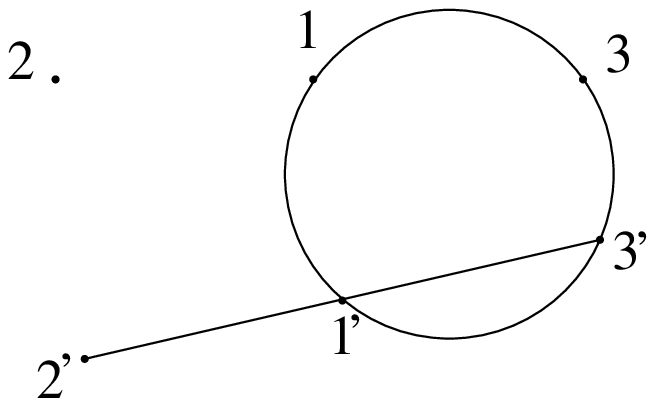, height=2cm}(b)
\caption{{\rm (a): Original theorem;\  (b): Two constraints removed.}}
\end{figure}

In \cite{li04} one hypothesis of the theorem was removed and a geometric factorization
of the conclusion was obtained. In \cite{li06} two hypotheses were removed and a geometric completion
was discovered. The computing was 3-termed in \cite{li04} and 2-termed in \cite{li06}. Both required
Clifford factorization, \cite{li04} further used circular transform.
Below we present a 1-termed computing without Clifford factorization nor transform.

The configuration in Figure 2b is constructed as follows:

\vskip -.23cm
Free points:
$\1,\2,\3,\1',\2'$. \\
Intersection: $\3'=\1'\2'\cap \1\3\1'$.\\
Conclusion expression: $[\e\2\2'\e\3\3']$.
\be
\hskip -.1cm
\ba{lll}
{[}\e\2\2'\e\3\3']&
=& 2^{-1}\,[\e\2\2'\e\3\{(\e\wedge \2')\vee_{\1'}(\1\wedge \3)\}\\
&&\hfill
\1'\{(\e\wedge \2')\vee_{\1'}(\1\wedge \3)\}]\\
&=& \underbrace{2^{-1}\,[\e\1\3\1'][\e\3\1'\2']}[\e\2\2'\e\3\1\1'\2']  \\
&=&  \underbrace{2\,(\e\cdot \2')}[\2'\2\e\3\1\1'].
\ea
\label{fnext}
\ee
BREEFS in (\ref{fnext}):

(1) The first meet product has two neighbors in the long geometric product: $\3$ and $\1'$. Expansion
by separating $\1,\3$ would simply
replace the meet product by $\1$. Separating $\e, \2'$ would produce two terms
and a Clifford factorization (\ref{cga:miquel4}) must be used to get a monomial result.

(2) Similarly,
the expansion of the second meet product has to separate $\e,\2'$, as $\e$ is a neighbor
by shift symmetry. The last step uses null symmetry
(\ref{cga:nullsymmetry}) and null expansion (\ref{cga:nullexpansion}):
\[\ba{ll}
[\e\2\2'\e\3\1\1'\2']\ =-[\2'\e\2\2'\e\3\1\1']&=[\2'\e\2'\2\e\3\1\1']\\
&\hfill
=2\,(\e\cdot \2')[\2'\2\e\3\1\1'].
\ea
\]

The result is not what we expected. The two brackets $[\e\1\2\3]$ and $[\1\2\1'\2']$ representing the two
removed hypotheses are not in the final result. They can be produced by
rational expansion (\ref{cga:miquel4}):
\[
\frac{1}{2}[\e\3\1\1'\2'\2]=\frac{\1\cdot \3[\e\2\3\2'][\1\2\1'\2']-\2\cdot \2'[\e\1\2\3][\1\3\1'\2']}
{[\1\2\3\2']}.
\]

{\bf Example 3.} [See \cite{li04}, Example 3]\ Let $\bf 1',2',3'$ be points on
sides $\bf 23, 13, 12$ of triangle $\bf 123$ respectively.
Then circles
$\bf 12'3'$, $\bf 1'23'$ and $\bf 1'2'3$ meet at a common point $\4$.

\begin{figure}[htbp]
\centerline{
\epsfig{file=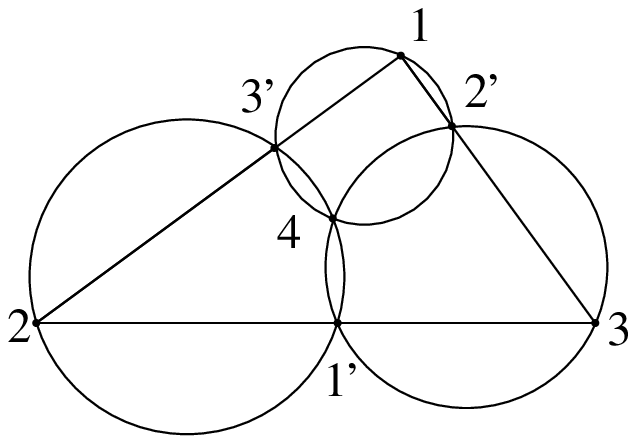, height=2.5cm}
\hskip 1cm
\epsfig{file=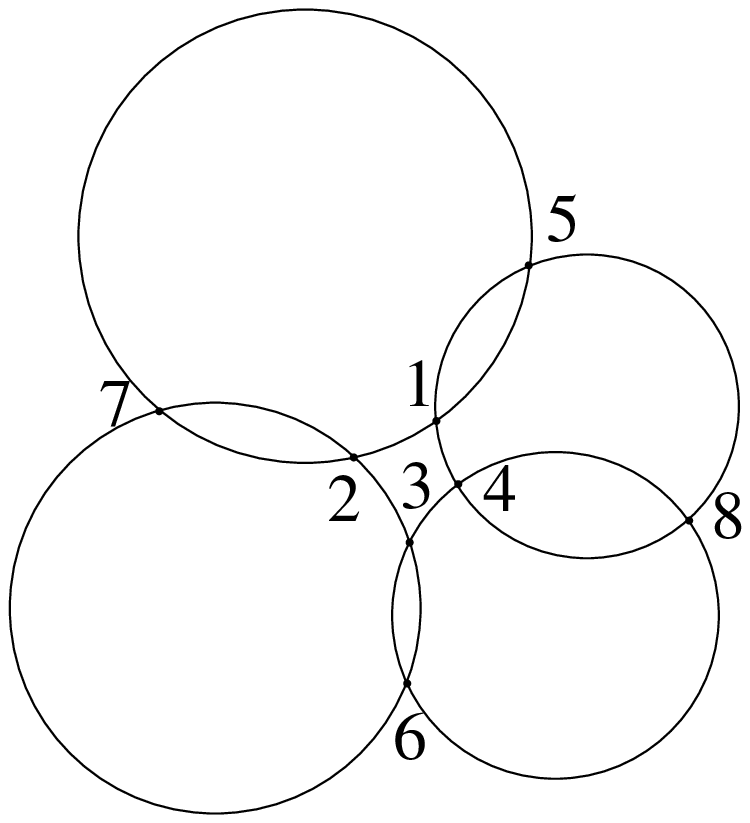, height=3.14cm}
}
\caption{Left: Example 3;\hskip .5cm Right: Example 4.}
\end{figure}

In \cite{li04} a collinearity constraint was removed in order to explore the
dependency of the conclusion upon the constraint. Below we remove
all the hypotheses (the three collinearity constraints) to explore the equivalence of the conclusion
with the hypotheses.

Conclusion: three circles $\1\2'\3'$, $\1'\2\3'$, $\1'\2'\3$ concur, i.e., the intersection
$\4=\1\2'\3'\cap \1'\2\3'$ is on circle $\1'\2'\3$:
\be
\hskip -.1cm
\ba{ll}
& [\1'\2'\3\4]\\
=& {[}\1'(\1\2'\3'\cap \1'\2\3')\2'\3]\\
=& \hskip -.13cm 2^{-1}[\1'\{(\1\wedge \2')\vee_{\3'}(\1'\wedge \2)\}
\3'\{(\1\wedge \2')\vee_{\3'}(\1'\wedge \2)\}\2'\3]\\
=& \hskip -.13cm\underbrace{2^{-1}[\1\1'\2'\3'][\2\1'\2'\3']}[\1'\2\3'\1\2'\3].\\
\ea
\label{cga:ex3new}
\ee
In conclusion expression $[\1'\2'\3\4]$,
$\1',\2',\3$ are antisymmetric. Because $\3$ is
irrelevant to $\4=\1\2'\3'\cap \1'\2\3'$, only when $\4$ is between $\1',\2'$
can the BREEFS principle take effect. Then similar to Example 2,
neighborhood consideration leads to unique monomial
expansions of the meet products.

Discarded hypotheses: three ``circles" $\1'\2\3$, $\1\2'\3$, $\1\2\3'$ concur (at ``point" $\e$).
This can be represented by the incidence of the intersection $\1'\2\3\cap \1\2'\3$ and
circle $\1\2\3'$. Simply by interchanging the primes over the same letters, we get from
(\ref{cga:ex3new}) the same effective part of the discarded hypotheses:
\be
[\1\2'\3\1'\2\3']=-[\1'\2\3'\1\2'\3]. \label{cga:ex33}
\ee
(\ref{cga:ex33})
discloses the intrinsic equivalence between the conclusion and the hypotheses.

{\bf Example 4.} [Miquel's 4-circle Theorem] Four circles intersect sequentially at
pairs of points
$(\1,\5)$, $(\2,\7)$, $(\3,\6)$ and $(\4,\8)$. If $\1,\2,\3,\4$ are cocircular then so are
$\5,\6,\7,\8$.

This is a typical theorem whose analytic proof using coordinates is difficult although straightforward.
In \cite{li02} a 5-termed NBA proof was found. Below we present a 1-termed proof.

Free points: $\bf 1,2,3,4,5,6$.\\
Intersections:
$
\7=\1\2\5\cap \2\3\6,\ \ \
\8=\1\4\5\cap \3\4\6.
$
\be
\hskip -.44cm
\ba{lrl}
&{[}\5\6\7\8] =& \hskip -.12cm-[\5\7\6\8] \\
&=& \hskip -.12cm-2^{-2}[\5\{(\1\wedge \5)\vee_\2(\3\wedge \6)\}\2
\{(\1\wedge \5)\vee_\2(\3\wedge \6)\} \\
&&\hskip .735cm \6\{(\1\wedge \5)\vee_\4(\3\wedge \6)\}\4
\{(\1\wedge \5)\vee_\4(\3\wedge \6)\}] \\
&=&\hskip -.12cm \bigstrut \underbrace{2^{-2}[\1\2\5\6][\1\4\5\6][\2\3\5\6][\3\4\5\6]}[\5\1\2\3\6\3\4\1]\\
&=&\hskip -.12cm \bigstrut
\underbrace{-2^2\,(\1\cdot \5)(\3\cdot \6)}[\1\2\3\4].
\ea
\label{miq4}
\ee
Similar to Example 3, neighbors of $\7,\8$ in bracket $[\5\6\7\8]$
should be changed to $\5,\6$ to make the best of the BREEFS principle. The rest are two simple
null expansions leading to the conclusion $[\1\2\3\4]$. If using the original $[\5\6\7\8]$, then
the two meet products
\[
\{(\1\wedge \5)\vee_\2(\3\wedge \6)\}\{(\1\wedge \5)\vee_\4(\3\wedge \6)\}
\]
are neighbors in the long geometric product. They
share two pairs of points: $(\1,\5)$ and $(\3,\6)$, and should be
expanded by separating the same pair. The result is 2-termed, and
Clifford factorization (\ref{cga:miquel4}) must be used to return to one term.

Summary of ``BREEFS" in NBA:

(1) To get factored and shortest result,
choose suitable algebraic representations, do eliminations and expansions in
long geometric products according to neighborhood consideration.
Immediate neighbors have topmost priority, then
next immediate neighbors, and so on. (2) If neighbors of a meet product
can be altered, then rearrange neighbors
by relevance to the meet product. (3) If an outer product is common to neighboring meet products,
expand the meet products by splitting common outer product.

Below we discuss Clifford factorizations (\ref{cga:miquel4}) to (\ref{cga:miquelformula}) which are
inverse of rational Clifford expansions.
They are the most important devices to reduce the number of terms of bracket polynomials. In particular,
the first formula of (\ref{cga:miquel4}) is the MOST frequently used Clifford factorization.
Its proof is straightforward
substitution of the Cramer's rule of $\4$ (or $\1$) with respect to basis $\2,\3,\5,\6$:
\be
[\2\3\5\6]\4=-[\3\4\5\6]\2+[\2\4\5\6]\3+[\2\3\4\6]\5-[\2\3\4\5]\6, \label{cga:4}
\ee
into bracket $[\1\2\3\4\5\6]$ and then null expansion (\ref{cga:nullexpansion}). All
other rational Clifford expansion formulas are easily proved in this way. However, the proof
does not provide any explanation as to why from the left to right of (\ref{cga:miquel4}) one
should add the factor $[\2\3\5\6]$, and why the right side should
contain any of the factors on the left.

(\ref{cga:miquel4}) from left to right can be understood
from the aspect of {\it pseudodivision}. In Example 2,
$[\1\2\3\4\5\6]$ in (\ref{cga:miquel4}) needs to be pseudodivided by
$[\3\4\5\6]$ and $[\1\2\3\6]$ to get the quotients. If $\4$ is the leading
vector variable, then since divisor $[\3\4\5\6]$ is linear with respect to $\4$, we can substitute
the expression of $\4$ by $\3,\5,\6$ into the dividend to get a remainder.
Similarly, if $\1$ is the leading vector variable and $[\1\2\3\6]$ is the divisor, we can
substitute the expression of $\1$ by $\2,\3,\6$ into the dividend to get another remainder.
The unique common pseudocoefficient for both pseudodivisions is $[\2\3\5\6]$. Thus the quotients
can be obtained by substituting (\ref{cga:4}) into the dividend.

To understand (\ref{cga:miquel4}) from right to left is much more difficult. However, once
a factor say $[\2\3\5\6]$ is discovered from the right, e.g., by polynomial
factorization in homogeneous coordinates, then the other factor $f$ can be easily fixed. Since
$f$ is a multilinear function of its vector variables, a procedure similar
to multilinear Cayley factorization \cite{white} leads to the result that $f$ equals $[\1\2\3\4\5\6]$
up to a constant factor.

\section{Elegance of analytic proof}
\vskip .3cm

{\bf Example 5.} [Miquel's 5-circle Theorem] Let there be a five star with vertices
$\bf 1', 2', 3', \4', \5'$, and armpit points $\bf 1, 2, 3, 4$, $\5$. The circumcircles of the five triangular
wedges meet sequentially at shoulder points $\bf 1'',2'',3'',4'',5''$ respectively.
Then the shoulder points are cocircular.
\begin{figure}[htbp]
\centerline{\epsfig{file=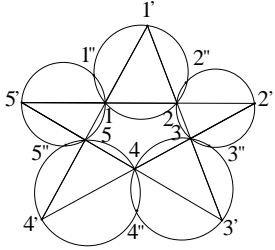, height=3.3cm}}
\caption{Miquel's 5-circle Theorem.}
\end{figure}

\vskip -.4cm
There is another saying is that this theorem is
due to W. K. Clifford, the inventor of Clifford algebra. Its
analytic proof is straightforward but complicated symbolic computing.
In this section we present a 3-termed beautiful proof with no short of
elegance than traditional synthetic one.

We use the same construction of the configuration as in \cite{li05}: first the armpit points
$\1,\2,\3,\4,\5$ as free points; then the vertices as intersections of lines:
\[\ba{lll}
\1'=\2\3\cap \5\1,&
\2'=\1\2\cap \3\4,&
\3'=\2\3\cap \4\5,\\
\4'=\3\4\cap \5\1,&
\5'=\4\5\cap \1\2;
\ea\]
finally the shoulder points as intersections of circles:
\[\ba{lll}
\1''=\1\1'\2\cap \5\5'\1, &
\2''=\2\2'\3\cap \1\1'\2, &
\3''=\3\3'\4\cap \2\2'\3, \\
\4''=\4\4'\5\cap \3\3'\4, &
\5''=\5\5'\1\cap \4\4'\5.
\ea
\]
By symmetry we only need to prove $[\1''\2''\3''\4'']=0$.
\be\ba{ll}
& [\1''\2''\3''\4'']\\
=& -2^{-4}
[\1\{(\1'\wedge \2)\vee_\1(\5\wedge \5')\}\{(\2'\wedge \3)\vee_\2(\1\wedge \1')\} \\
&\hskip .85cm \2\{(\2'\wedge \3)\vee_\2(\1\wedge \1')\}\{(\3'\wedge \4)\vee_\3(\2\wedge \2')\}\\
&\hskip .85cm \3\{(\3'\wedge \4)\vee_\3(\2\wedge \2')\}\{(\4'\wedge \5)\vee_\4(\3\wedge \3')\}\\
&\hskip .85cm \4\{(\4'\wedge \5)\vee_\4(\3\wedge \3')\}\{(\1'\wedge \2)\vee_\1(\5\wedge \5')\}].
\ea
\label{cga:miquel5:0}
\ee

Consider the expansion of the first meet product. Its immediate neighbors provide no hint.
A next immediate neighbor is $\2$, suggesting the separation of $\1',\2$. Similarly, for the second meet
product, its next immediate neighbor $\1$ suggests
the expansion by separating $\1,\1'$.
The two expansions result in three terms from the first two meet products:
\[\ba{ll}
& \1\{(\1'\wedge \2)\vee_\1(\5\wedge \5')\}\{(\2'\wedge \3)\vee_\2(\1\wedge \1')\}\2\\
=& -[\1\5\1'\5'][\2\3\1'\2'] \1\2\1\2\bigstrut
+[\1\2\3\2'][\1\5\1'\5']\1\2\1'\2\\
&\hfill -[\1\2\5\5'][\2\3\1'\2']\1\1'\1\2  \\
\ea\]
\[\ba{ll}
=& 2\,(-\1\cdot \2[\1\5\1'\5'][\2\3\1'\2']\bigstrut
+\2\cdot \1'[\1\2\3\2'][\1\5\1'\5']\\
&\hfill
-\1\cdot \1'[\1\2\5\5'][\2\3\1'\2'])\1\2\\
=& g_{12}\1\2.
\ea
\]

\vskip -.3cm
Suddenly it appears that the first two meet products are simply removed from (\ref{cga:miquel5:0}):
\be\ba{ll}
& [\1''\2''\3''\4'']\\
=& \ds \hskip -.2cm\bigstrut
\underbrace{-2^{-4} g_{12}}[\1\2\{(\2'\wedge \3)\vee_\2(\1\wedge \1')\}\{(\3'\wedge \4)\vee_\3(\2\wedge \2')\}\\
&\hskip 1.33cm \3\{(\3'\wedge \4)\vee_\3(\2\wedge \2')\}\{(\4'\wedge \5)\vee_\4(\3\wedge \3')\}\\
&\hskip 1.33cm\4\{(\4'\wedge \5)\vee_\4(\3\wedge \3')\}\{(\1'\wedge \2)\vee_\1(\5\wedge \5')\}].
\ea
\label{cga:miquel5:1}
\ee
Although irrelevant to the proof,
It is interesting to note that after eliminating $\1',\2',\5'$ from its expression, the
intermediate factor $g_{12}$ equals
\[
-2\,(\1\cdot \2)(\1\cdot \5)(\2\cdot \3)[\e\1\2\3]^2[\e\1\2\5]^3[\e\1\4\5][\e\2\3\4][\e\3\4\5]
S_{\bf 13524},
\]
where
\[\ba{lll}
S_{\bf 13524}&=& (\e\cdot \1)(\e\cdot \4)[\e\2\3\5]
-(\e\cdot \1)(\e\cdot \5)[\e\2\3\4]\\
&&\hfill
+(\e\cdot \2)(\e\cdot \5)[\e\1\3\4]
\ea
\]
is twice the signed area of pentagon $\1\3\5\2\4$.

Because of the symmetry in the geometric constructions, in (\ref{cga:miquel5:1})
the expansion of the two meet products
between $\2,\3$ results in an intermediate factor $g_{23}$, and the expansion of the two meet products
between $\3,\4$ results in an intermediate factor $g_{34}$. Of course by direct
computing we obtain the same result.
Now (\ref{cga:miquel5:1}) is changed into
\be
h=\underbrace{g_{23}g_{34}}[\1\2\3\4\{(\4'\wedge \5)\vee_\4(\3\wedge \3')\}
\{(\1'\wedge \2)\vee_\1(\5\wedge \5')\}].
\label{cga:mq5}
\ee

In (\ref{cga:mq5}), immediate neighbors of the meet products
suggest two ways of expanding them simultaneously, each resulting
in three terms: either separate $\4',\5$ and $\5,\5'$, or separate
$\3, \3'$ and $\1', \2$. The latter expansion has three terms
because
\[
[\1\2\3\4\3\2]=2\,(\3\cdot \4)[\1\2\3\2]=0.
\]
To choose between the two options, consider the next immediate neighbors.
The first meet product has next immediate neighbors $\3, \1$, and the second
has $\2, \4$. They suggest the separation of
$\3, \3'$ and $\1',\2$ respectively.
The benefit is immediate null expansion of long brackets:
\be
\hskip -.42cm
\ba{lrl}
& h =&\hskip -.17cm -[\1\2\5\5'][\4\5\3'\4'][\1\2\3\4\3\1']+[\1\5\1'\5'][\3\4\5\4'][\1\2\3\4\3'\2]\\
&& \hfill +\, [\1\2\5\5'][\3\4\5\4'][\1\2\3\4\3'\1']\\
&=& \bigstrut
\hskip -.17cm-2\,(\3\cdot \4)[\1\2\3\1'][\1\2\5\5'][\4\5\3'\4']-2\,(\1\cdot \2)[\1\5\1'\5']\times\\
&&\hfill \times
[\2\3\4\3'][\3\4\5\4']+[\1\2\5\5'][\3\4\5\4'][\1\2\3\4\3'\1'].
\ea
\label{cga:miquel5:2}
\ee

The first round of elimination finishes here. In the second round, all remaining
constrained points are eliminated:
\[\ba{lll}
{[}\1\2\3\1'] &\stackrel{\1'}{=}& -\2\cdot \3[\e\1\2\3][\e\1\2\5][\e\1\3\5], \\
{[}\1\2\5\5'] &\stackrel{\5'}{=}& -\1\cdot \2[\e\1\2\5][\e\1\4\5][\e\2\4\5], \\
{[}\2\3\4\3'] &\stackrel{\3'}{=}& -\2\cdot \3[\e\2\3\4][\e\2\4\5][\e\3\4\5], \\
{[}\3\4\5\4'] &\stackrel{\4'}{=}& -\3\cdot \4[\e\1\3\5][\e\1\4\5][\e\3\4\5], \\
{[}\1\5\1'\5'] &\stackrel{\1'}{=}& \ \ \,\1\cdot \5[\e\1\2\3][\e\1\5\5'][\e\2\3\5], \\
{[}\e\1\5\5'] &\stackrel{\5'}{=}& \ds -2^{-1} [\e\1\2\e\4\5][\e\1\2\5][\e\1\4\5],\\
{[}\4\5\3'\4'] &\stackrel{\3'}{=}& \ \ \,\4\cdot \5[\e\2\3\4][\e\2\3\5][\e\4\5\4'], \\
\ea\]
\[\ba{lll}
{[}\e\4\5\4'] &\stackrel{\4'}{=}& \ds -2^{-1} [\e\3\4\e\5\1][\e\1\4\5][\e\3\4\5].
\ea
\]
Then
\be\ba{ll}
h =& \underbrace{(\1\cdot \2)(\3\cdot \4)[\e\1\2\5][\e\1\3\5][\e\1\4\5]^2[\e\2\4\5][\e\3\4\5]}\\
&\hskip .5cm \{
\2\cdot \3[\e\1\2\3][\e\2\3\4][\e\2\3\5](\4\cdot \5[\e\1\2\5][\e\3\4\e\5\1]\\\
&\hfill +\1\cdot \5[\e\3\4\5][\e\1\2\e\4\5])+[\1\2\3\4\3'\1']\}.
\ea
\label{cga:miquel5:3}
\ee

(\ref{cga:miquel5:3}) contains a matching with (\ref{cga:miquelformula}) for Clifford factorization,
$\1\2\3\4\5\6$ in the first formula for $\e\5\4\3\1\2$ here:
\be\ba{r}
\4\cdot \5[\e\1\2\5][\e\3\4\e\5\1]+
\1\cdot \5[\e\3\4\5][\e\1\2\e\4\5]\\
\ds
=-\frac{1}{2}[\e\5\4\3\e\5\1\2][\e\1\4\5].
\ea
\ee
Then (\ref{cga:miquel5:3}), and consequently the conclusion, is reduced to
\be
[\1\2\3\4\3'\1']=2^{-1}(\2\cdot \3)[\e\1\2\3][\e\1\4\5][\e\2\3\4][\e\2\3\5][\e\5\4\3\e\5\1\2].
\label{cga:miquel5:4}
\ee

Below we prove (\ref{cga:miquel5:4}).
\[\ba{ll}
& {[}\1\2\3\4\3'\1'] \\
=& 2^{-2}[\1\2\3\4\{(\2\wedge \3)\vee_\e(\4\wedge \5)\}\e\bigstrut
\{(\2\wedge \3)\vee_\e(\4\wedge \5)\}\\
& \hskip 1.34cm  \{(\2\wedge \3)\vee_\e(\1\wedge \5)\}\e
\{(\2\wedge \3)\vee_\e(\1\wedge \5)\}].
\ea\]
By neighborhood consideration, the first meet product is expanded by
separating $\4,\5$, the last meet product is expanded by
separating $\1, \5$:
\be\hskip -.5cm
\ba{l}
{[}\1\2\3\4\3'\1'] =2^{-2}[\e\1\2\3][\e\2\3\4][\1\2\3\4\5\e\{(\2\wedge \3)\vee_\e\\
\hskip 3cm(\4\wedge \5)\}
\{(\2\wedge \3)\vee_\e(\1\wedge \5)\}\e\5].
\ea
\label{cga:miquel5:meet}
\ee

The meet products in (\ref{cga:miquel5:meet}) form the double-line type in Cayley expansion theory
\cite{li-wu}:
\be\ba{ll}
& \e\{(\2\wedge \3)\vee_\e(\4\wedge \5)\}\{(\2\wedge \3)\vee_\e(\1\wedge \5)\}\e\\
=& \e\{\{(\2\wedge \3)\vee_\e(\4\wedge \5)\}\wedge \{(\2\wedge \3)\vee_\e(\1\wedge \5)\}\}\e\\
=& \{(\2\wedge \3)\vee_\e(\4\wedge \5)\vee_\e(\1\wedge \5)\}\,\e(\2\wedge \3)\e\\
=& -[\e\1\4\5][\e\2\3\5]\e\2\3\e.
\label{cga:miquel5:lem2}
\ea\ee
Then (\ref{cga:miquel5:4}) is equivalent to
\be
[\1\2\3\4\5\e\2\3\e\5]=-2\,(\2\cdot \3)[\e\5\4\3\e\5\1\2].
\label{cga:miquel5:6}
\ee
It must be pointed out that without using Cayley expansion we can simply expand the meet products
in (\ref{cga:miquel5:meet}) by separating $\2,\3$ simultaneously. The result is 2-termed which can
be contracted to one term using Grassmann-Pl\"ucker syzygy (\ref{cga:gp}). The proof still remains 3-termed.

Now Miquel's 5-circle theorem is equivalent to algebraic identity (\ref{cga:miquel5:6}).
Use null symmetry and shift symmetry to rearrange the sequence so that two
$\3$'s are separated by three vectors,
then use the trigonometric quartet expansion and factorization (\ref{quartet}), and
null expansion (\ref{cga:nullexpansion}) to get
\be\ba{ll}
& {[}\1\2\3\4\5\e\2\3\e\5] \\
=& -[\3\4\5\e\3\2\e\5\1\2] \\
=& -2\,(\langle \3\4\5\e\rangle \,[\3\2\e\5\1\2]+[\3\4\5\e]\langle \3\2\e\5\1\2\rangle ) \\
=& -2^2\,(\2\cdot \3)(\langle \e\5\4\3\rangle \,[\e\5\1\2]+[\e\5\4\3]\,\langle \e\5\1\2\rangle)\\
=& -2\,(\2\cdot \3)[\e\5\4\3\e\5\1\2].
\ea\ee

This finishes the proof of the theorem.
Summing up, with $(\vee_{\bf i})$ denoting a reduced meet product with respect to $\bf i$,
the proving of the whole theorem proceeded as follows:

\[\hskip -.13cm
\ba{ll}
{[}\1''\2''\3''\4'']&\hskip -.13cm
\longrightarrow
[\1(\vee_\1)(\vee_\2)\2(\vee_\2)(\vee_\3)\3(\vee_\3)(\vee_\4)\4(\vee_\4)(\vee_\1)]\\
&\hskip -.13cm\longrightarrow {[}\1\2(\vee_\2)(\vee_\3)\3(\vee_\3)(\vee_\4)\4(\vee_\4)(\vee_\1)]\\
&\hskip -.13cm\longrightarrow {[}\1\2\3(\vee_\3)(\vee_\4)\4(\vee_\4)(\vee_\1)]\\
&\hskip -.13cm\longrightarrow {[}\1\2\3\4(\vee_\4)(\vee_\1)]\\
&\hskip -.13cm\longrightarrow \hbox{3-termed } h \hbox{ in (\ref{cga:miquel5:2})} \\
&\hskip -.13cm\longrightarrow \hbox{3-termed } h \hbox{ in (\ref{cga:miquel5:3})} \\
&\hskip -.13cm\longrightarrow \hbox{2-termed (\ref{cga:miquel5:4}): } [\1\2\3\4\3'\1'],\ [\e\5\3\4\e\5\1\2],
\ea
\]
succeeded by
\[\ba{ll}
{[}\1\2\3\4\3'\1']&
\longrightarrow
{[}\1\2\3\4(\vee_\e)\e(\vee_\e)(\vee_\e)\e(\vee_\e)]\\
&\longrightarrow {[}\1\2\3\4\5\e(\vee_\e)(\vee_\e)\e\5]\\
&\longrightarrow {[}\1\2\3\4\5\e\2\3\e\5]\\
&\longrightarrow {[}\e\5\3\4\e\5\1\2].
\ea\]

\end{document}